\documentclass{amsart}

\newtheorem{theorem}{Theorem}[section]
\newtheorem{lemma}[theorem]{Lemma}
\newtheorem{proposition}[theorem]{Proposition}
\newtheorem{corollary}[theorem]{Corollary}

\theoremstyle{definition}
\newtheorem{definition}[theorem]{Definition}

\newtheorem{example and remarks}[theorem]{Example and Remarks}
\newtheorem{example and remark}[theorem]{Example and Remark}
\newtheorem{notation}[theorem]{Notation}
\newtheorem{remarks}[theorem]{Remarks}
\newtheorem{notations}[theorem]{Notations}
\newtheorem{remark}[theorem]{Remark}

\numberwithin{equation}{section}

\begin{document}

\title[ON THE LOCAL COHOMOLOGY MODULES DEFINED BY A PAIR OF IDEALS ]{ ON THE LOCAL COHOMOLOGY MODULES DEFINED BY A PAIR OF IDEALS AND SERRE SUBCATEGORIES }

\author{KH.  AHMADI-AMOLI}
\address{Department of Mathematics, Payame Noor university, Tehran, 19585-83133, Iran}
\curraddr{Department of Mathematics, Payame Noor university, Tehran, 19585-83133, Iran}
\email{khahmadi@pnu.ac.ir}
\thanks{}

\author{M.Y.  SADEGHI}
\address{Department of Mathematics, Payame Noor university, Tehran, 19585-83133, Iran}
\email{sadeghi@phd.ac.ir}
\thanks{}


\date{}

\dedicatory{ }


\begin{abstract}
This paper is concerned about the relation between local cohomology modules defined by a pair of ideals and Serre classes of $R$-modules, as a generalization of results of $~$J. Azami, R. Naghipour and B. Vakili (2009) and M. Asgharzadeh and M.Tousi (2010). Let $R$ be
a commutative Noetherian ring, $I$ , $J$ be two ideals of $R$ and $M$ be an $R$-module. Let  $\mathfrak{a}\in\tilde{W}(I,J)$ and $t\in\mathbb{N}_0$ be such that Ext$^t_R(R/\mathfrak{a},M)\in\mathcal{S}$ and Ext$^j_R(R/\mathfrak{a},H^i_{I,J}(M))\in\mathcal{S}$ for all $i<t$ and all $j\geq0$. Then for any submodule $N$ of $H^t_{I,J}(M)$ such that Ext$^1_R(R/\mathfrak{a},N)\in\mathcal{S}$, we obtain Hom$_R(R/\mathfrak{a},H^t_{I,J}(M)/N)$
$\in$ $\mathcal{S}$.
\end{abstract}

\maketitle

\section{Introduction}

Throughout this paper, $R$ is denoted a commutative Noetherian ring, $I$ , $J$ are denoted two ideals of $R$, and $M$ is denoted an arbitrary $R$-module. By $\mathbb{N}_0$, we shall mean the set of non-negative integers. For basic results, notations and terminologies not given in this paper, the reader is referred to [7] and [22], if necessary.\\
As a generalization of the usual local cohomology modules, Takahashi, Yoshino and Yoshizawa [22], introduce the local cohomology modules with respect to a pair of ideals $(I,J)$. To be more precise, let   $W(I,J)$ = $\{~\mathfrak{p}$ $\in$ Spec$(R)$ $\mid$ $I^n\subseteq \mathfrak{p}+J$ for some positive integer $n\}$ and $\tilde{W}(I,J)$ denote the set of ideals $\mathfrak{a}$ of $R$ such that $I^n\subseteq \mathfrak{a}+J$ for some integer $n$. In general, $W(I,J)$ is closed under specialization, but not necessarily a closed subset of Spec$(R)$. For an $R$-module $M$, we consider the $(I,J)$-torsion submodule $\Gamma_{I,J}(M)$ of $M$ which consists of all elements $x$ of $M$ with Supp$(Rx)$ in $W(I,J)$. Furthermore, for an integer $i$, the local cohomology functor $H^i_{I,J}$ with respect to $(I,J)$ is defined to be the $i$-th right derived functor of $\Gamma_{I,J}$. Also $H^i_{I,J}(M)$ is called the $i$-th local cohomology module of $M$ with respect to $(I,J)$. If $J$ = 0, then $H^i_{I,J}$ coincides with the ordinary local cohomology functor $H^i_I$ with the support in the closed subset $V(I)$.

Recently, some authors approached the study of properties of these extended modules, see for example [9], [10], [19] and [23].

It is well known that an important problem in commutative algebra is to determine when the $R$-module Hom$_R(R/I,H^i_I(M))$ is finite. Grothendieck in [14] conjectured the following:

If $R$ is a Noetherian ring, then for any ideal  $I$ of $R$ and any finite $R$-module $M$, the modules Hom$_R(R/I,H_I^i(M))$
are finite for all $i\geq0$.

In [15], Hartshorne gave a counterexample to Grothendieck$^,s$ conjecture and he defined the concept of $I$-cofinite modules to generalize the conjecture. In [6], Brodmann and Lashgari showed that if, for a finite $R$-module $M$ and an integer $t$, the local cohomology modules $H^0(M)$ , $H^1(M)$ , $\cdots$, $H^{t-1}(M)$ are finite, then $R$-module Hom$_R(R/I,H^t_I(M))$ is finite and so Ass$(H^t_I(M)/N)$ is a finite set for any finite submodule $N$ of $H^t_I(M)$. A refinement of this result for $I$-minimax $R$-modules is as follows, see [4].

\begin{theorem} \label{1.1}
Let $M$ be an $I$-minimax $R$-module and $t$ be a non-negative integer such that $H^i_I(M)$ is $I$-minimax for all $i<t$. Then for any $I$-minimax submodule $N$ of $H^t_I(M)$, the $R$-module Hom$_R(R/I,H^t_I(M))$ is $I$-minimax.
\end{theorem}

Also authors in [1] and [3] studied local cohomology modules by means of Serre subcategories. As a consequence, for an arbitrary Serre subcategory $\mathcal{S}$, authors in [3] showed the following result.

\begin{theorem} \label{1.2}
Let $s \in \mathbb{N}_0$ be such that Ext$^s_R(R/I,M)\in\mathcal{S}$ and Ext$^j_R(R/I,H^i_I(M))\\\in\mathcal{S}$ for all $i<s$ and all $j\geq0$. Let $N$ be a submodule of $H^s_I(M)$ such that Ext$^1_R(R/I,N)\in\mathcal{S}$. Then Hom$_R(R/I,H^s_I(M)/N)\in\mathcal{S}$.
\end{theorem}

The aim of the present paper is to generalize the concept of $I$-cominimax $R$-module, introduced in [4], to an arbitrary Serre subcategory $\mathcal{S}$, to verify situations in which the $R$-module Hom$_R(R/I,H^i_{I,J}(M))$ belongs to $\mathcal{S}$. To approach it, we use the methods of [3] and [4]. Our paper consists of four sections as follows.

In Section 2, by using the concept of $(I,J)$-relative Goldie dimension, we introduce the $(I,J)$-minimax $R$-modules and we study some properties of them, see Proposition {2.7}.

In Section 3, for an arbitrary Serre subcategory $\mathcal{S}$, we defined $(\mathcal{S},I,J)$-cominimax $R$-modules. This concept of $R$-modules can be considered as a generalization of $I$-cofinite $R$-modules [14], $I$-cominimax $R$-modules [4], and $(I,J)$-cofinite $R$-modules [23]. Also, as a main result of our paper, we prove the following. (See Theorem {3.4}).

\begin{theorem} \label{1.3}
Let  $\mathfrak{a}\in\tilde{W}(I,J)$. Let $t\in\mathbb{N}_0$ be such that Ext$^t_R(R/\mathfrak{a},M)$ $\in$ $\mathcal{S}$ and Ext$^j_R(R/\mathfrak{a},H^i_{I,J}(M))\in\mathcal{S}$ for all $i<t$ and all $j\geq0$. Then for any submodule $N$ of $H^t_{I,J}(M)$ such that Ext$^1_R(R/\mathfrak{a},N)\in\mathcal{S}$, we have Hom$_R(R/\mathfrak{a},H^t_{I,J}(M)/N)\in\mathcal{S}$.
\end{theorem}

One can see, by replacing various Serre classes with $\mathcal{S}$ and using Theorem {1.3}, the main results of [{2}, Theorem 1.2] , [{3}, Theorem 2.2] , [{4}, Theorem 4.2] , [{5}, Lemma 2.2] , [6] , [{12}, Corollary 2.7] , [16] , [{17}, Corollary 2.3], and [{23}, Theorem 3.2] are obtained.(See Theorem {3.14} and Proposition {3.15}).

At last, in Section 4, as an application of results of the previous sections, we give the following consequence about finiteness of associated primes of local cohomology modules.(See Proposition {4.1} and Corollary {4.2})

\begin{proposition} \label {1.4}
Let $t\in\mathbb{N}_0$ be such that Ext$^t_R(R/I,M)\in\mathcal{S}_{I,J}$ and $H^i_{I,J}(M)\in\mathcal{C}(\mathcal{S}_{I,J},I,J)$ for all $i<t$. Let $N$ be a submodule of $H^t_{I,J}(M)$ such that Ext$^1_R(R/I,N)$ belongs to $\mathcal{S}_{I,J}$. If Supp$(H^t_{I,J}(M)/N)\subseteq V(I)$ then Gdim$H^t_{I,J}(M)/N<\infty$ and so $H^t_{I,J}(M)/N$ has finitely many associated primes; in particular, for $N$ = $JH^t_{I,J}(M)$.
\end{proposition}

\section{Serre Classes And $(I,J)$-Minimax Modules}

Recall that for an $R$-module $H$, the Goldie dimension of $H$ is defined as the cardinal of the set of indecomposable submodules of $E(H)$, which appear in a decomposition of $E(H)$ in to direct sum of indecomposable submodules. Therefore, $H$ is said to have finite Goldie dimension if $H$ does not contain an infinite direct sum of non-zero submodules, or equivalently the injective hull $E(H)$ of $H$ decomposes as a finite direct sum of indecomposable (injective) submodules. We shall use Gdim$~H$ to denote the Goldie dimension of $H$. For a prime ideal $\mathfrak{p}$, let $\mu^0(\mathfrak{p},H)$ denotes the $0$-th Bass number of $H$ with respect to prime ideal $\mathfrak{p}$.
It is known that $\mu^0(\mathfrak{p},H) > 0$  iff  $\mathfrak{p}\in Ass(H)$. It is clear by the definition of the Goldie dimension that Gdim$~H$ = $\sum_{\mathfrak{p}\in Ass(H)}$ $\mu^0(\mathfrak{p},H)$ = $ \sum_{\mathfrak{p}\in Spec(R)}$ $\mu^0(\mathfrak{p},H)$. Also, the $(I,J)$-relative Goldie dimension of $H$ is defined as Gdim$_{I,J}H$ {:}= $ \sum_{\mathfrak{p}\in W(I,J)}$ $\mu^0(\mathfrak{p},H)$.(See [{19}, Definition 3.1]). If $J$ = 0, then Gdim$_{I,J}H$ = Gdim$_{I}H$ (see[{11}, Definition 2.5]),  moreover if $I$ = 0, we obtain Gdim$_{I,J}H$ = Gdim$~H$. It is known that when $R$ is a Noetherian ring, an $R$-module $H$ is minimax if and only if any homomorphic image of $H$ has finite Goldie dimension (see [13], [25], or [26]). This motivates the definition of $(I,J)$-minimax modules.

\begin{definition} \label{2.1}
An $R$-module $M$ is said to be minimax with respect to $(I,J)$ or $(I,J)$-minimax if the $(I,J)$-relative Goldie  dimension of any quotient module of $M$ is finite.
\end{definition}

\begin{remarks} \label{2.2}
By Definition {2.1}, it is clear that
\begin{enumerate}
\item[(i)]
Gdim$_{I}M \leq$ Gdim$_{I,J}M \leq$ Gdim$~M.$
\end{enumerate}
This inequalities maybe strict (see [{19}, Definition 3.1]).
\begin{enumerate}
\item[(ii)]  For Noetherian ring $R$, an $R$-module $M$ is minimax iff for any $R$-module $N$ of $M$, Gdim$~M/N<\infty$. Therefore, by (i), in Noetherian case, the class of $I$-minimax $R$-modules contains the class of $(I,J)$-minimax $R$-modules and it contains the class of minimax $R$-modules.
\end{enumerate}
\end{remarks}

\begin{example and remarks} \label{2.3}
It is easy to see that
\begin{enumerate}
\item[(i)]  Every quotient of finite modules, Artinian modules, Matlis refelexive modules and linearly compact modules have finite $(I,J)$-relative Goldie dimension, and so all of them are $(I,J)$-minimax modules.

\item[(ii)]  If $I$ $=$ 0, then $W(I,J)$ = $Spec(R)$ = $V(I)$ and so an $R$-module $M$ is minimax iff is $(I,J)$-minimax iff is $I$-minimax.

\item[(iii)]  If $J$ $=$ 0, then $W(I,J)$ = $V(I)$, so that an $R$-module $M$ is $(I,J)$-minimax iff is $I$-minimax.

\item[(iv)]  let $M$ be an $I$-torsion module. Then, by [{11}, Lemma 2.6] and [{19}, Lemma 3.3], $M$ is minimax iff is $(I,J)$-minimax iff is $I$-minimax.

\item[(v)]  If $M$ is $(I,J)$-torsion module, then $M$ is minimax iff is $(I,J)$-minimax. (By the definition and [{19}, Lemma 3.3]). Specially, when $(0)$ $\in$ \~{W}$(I,J)$ or $(0)$ $\in$ W$(I,J)$.

\item[(vi)]  By [{22}, Corollary (1.8)(2)], the class of $(I,J)$-torsion is a Serre subcategory of $R$-modules. Therefore, by part(v), in this category, the concept of minimax modules coincides with the concept of $(I,J)$-minimax modules; specially, for the $(I,J)$-torsion module of $H_{I,J}^i(M)$ ($i\geq0$).

\item[(vii)] If $Min(M)$ $\subseteq$ $W(I,J)$ and Gdim$_{I,J}M<\infty$, then by [{19}, Lemma 3.3] and definition, Gdim$~M<\infty$, and so $|Ass(M)|<\infty$.
\end{enumerate}
\end{example and remarks}

The following proposition shows that the class of $(I,J)$-minimax $R$-modules is a Serre subcategory.

\begin{proposition} \label{2.4}
Let  $0\rightarrow {M^\prime}\rightarrow M$ $\rightarrow M^\prime$$^\prime$ $\rightarrow0$ be an exact sequence of $R$-modules. Then $M$ is $(I,J)$-minimax if and only if M$^\prime$ and M$^\prime$$^\prime$ are both $(I,J)$-minimax.
\begin{proof}
One can obtain the result, by replacing $W(I,J)$ with $V(I)$ and a modification of the proof of proposition {2.3} of [4].
\end{proof}
\end{proposition}

\begin{remark} \label{2.5}
Recall that a class $\mathcal{S}$ of $R$-modules is a $''$Serre  subcategory$''$ or $''$Serre class$''$ of the category of $R$-modules, when it is closed under taking submodules, quotients and extensions. For example, the following class of $R$-modules are Serre  subcategory.
\begin{enumerate}

\item[(a)] The class of Zero modules.

\item[(b)] The class of Noeterian modules.

\item[(c)] The class of Artinian modules.

\item[(d)] The class of $R$-modules with finite support.

\item[(e)] The class of all $R$-modules $M$ with $dim_R M\leq n$, where $n$ is a non-negative integer.

\item[(f)] The class of minimax modules and the class of $I$-cofinite minimax $R$-modules. ~(see [{18}, Corollary 4.4] )

\item[(g)] The class of $I$-minimax $R$-modules. ~(see [{4}, Proposition 2.3 ] )
\item[(h)] The class of $I$-torsion $R$-modules and the class of $(I,J)$-torsion $R$-modules. (see [{22}, Corollary 1.8])

\item[(i)] The class of $(I,J)$-minimax $R$-modules. ~(Proposition {2.4})
\end{enumerate}
\end{remark}

\begin{notations} \label{2.6}
In this paper, the following notations are used for the following Serre subcategories:
\begin{enumerate}
\item[] $\mathcal{''~S~''}$ for an arbitrary Serre class of $R$-modules;
\item[] $\mathcal{''~S}_0~''$ for the class of minimax $R$-modules;
\item[] $\mathcal{''~S}_I~''$ for the class of $I$-minimax $R$-modules;
\item[] $\mathcal{''~S}_{I,J}~''$ for the class of $(I,J)$-miniax $R$-modules.
\end{enumerate}
\end{notations}

Using the above notations and Remark {2.2}, we have $~\mathcal{S}$$_0$ $\subseteq\mathcal{S}$$_{I,J}$ $\subseteq\mathcal{S}$$_{I}$ .

Now, we exhibit some of the properties of  $\mathcal{S}_{I,J}$ .

\begin{proposition} \label{2.7}
Let $I$,$J$,$I^\prime$,$J^\prime$ be ideals of $R$ and $M$ be an $R$-module. Then
\begin{enumerate}

\item[(i)] $\mathcal{S}_{I,J}$ = $\mathcal{S}_{{\sqrt{I}},J}$ = $\mathcal{S}_{{\sqrt{I}},{\sqrt{J}}}$ = $\mathcal{S}_{I,{\sqrt{J}}}$ .

\item[(ii)]  If $I^n$$\subseteq{\sqrt{J}}$, for some $n\in\mathbb{N}$ (or equally, if $R/J$ is an $I$-torsion $R$-module), then $\mathcal{S}_{0}$ = $\mathcal{S}_{I,J}$ .

\item[(iii)]  If $I^n$ $\subseteq{\sqrt{I'}}$, for some $n\in\mathbb{N}$, then $\mathcal{S}_{I,J}$ $\subseteq$ $\mathcal{S}_{I^\prime,J}$ .

\item[(iv)]  If $J^n$ $\subseteq{\sqrt{J'}}$, for some $n\in\mathbb{N}$, then $\mathcal{S}_{I,J^\prime}$ $\subseteq$ $\mathcal{S}_{I,J}$ .

\item[(v)]  If $I^n$ $\subseteq{\sqrt{I'}}$, for some $n\in\mathbb{N}$ and $M$ is $(I^\prime,J)$-torsion, then $M \in$ $\mathcal{S}_{I,J}~$ iff $~M\in\mathcal{S}_{0}~$ iff $~M\in\mathcal{S}_{{I^\prime},J}$ .

\item[(vi)]  If $J^n$ $\subseteq{\sqrt{J'}}$, for some $n\in \mathbb{N}$ and $M$ is $(I,J)$-torsion, then $M \in$ $\mathcal{S}_{I,J}~$ iff $~M\in\mathcal{S}_{0}~$ iff $~M\in\mathcal{S}_{I,{J^\prime}}$ .
\end{enumerate}
\begin{proof}

All these statements follow easily from [{22}, Proposition 1.4 and 1.6] and Remark {2.3}. As an illustration, we just prove statement (iii).\\
 Let $H \in\mathcal{S}_{I,J}$. Since $I^n\subseteq\sqrt{I'}$, we have $W({\sqrt{I'}},J)\subseteq W(I,J)$, by [{22}, Proposition 1.6]. Now, since $H$ is $(I,J)$-minimax, the assertion follows from definition.
\end{proof}
\end{proposition}

\begin{lemma} \label{2.8}
(i) If $N\in\mathcal{S}$ and $M$ is a finite $R$-module, then for any submodule $H$ of $Ext_R^i(M,N)$ and $T$ of $Tor_R^i(M,N)$, we have  $Ext_R^i(M,N)/H\in\mathcal{S}$ and $Tor_R^i(M,N)/T\in\mathcal{S}$, for all $i\geq0$ .\\
(ii) For all $i\geq0$, we have $Ext_R^i(R/I,M)\in\mathcal{S}_{I,J}~$ iff $~Ext_R^i(R/I,M)\in\mathcal{S}_{0}~$ iff $~Ext_R^i(R/I,M)\in\mathcal{S}_{I}$ .
\begin{proof}
(i) The result follows from [{3}, Lemma 2.1].\\
(ii) Since, for all $i$, $Ext_R^i(R/I,M)$ and $Tor_i^R(R/I,M)$ are $(I,J)$-torsion $R$-modules, the assertion holds by Remark {2.3} (iv).
\end{proof}
\end{lemma}

The following proposition can be thought of as a generalization of Proposition {2.6} of [4], in case of $J$ = 0 and $\mathcal{S}$ = $\mathcal{S}_I$ .

\begin{proposition} \label{2.9}
Let $Min(M)\subseteq W(I,J)$. If $M~$$\in$$~\mathcal{S}$, then $H_{I,J}^i(M)~$$\in~$$\mathcal{S}$ for all $i\geq0$.
\begin{proof}
By hypothesis and [{21}, Corollary 1.7], $M$ is $(I,J)$-torsion $R$-module and so $H_{I,J}^0(M)~$=$~\Gamma_{I,J}$$(M)$~=~$M$. Therefore, $H_{I,J}^i(M)~$=~0 for all $i\geq1$, by [{22}, Corollary 1.13]. Thus the assertion holds.
\end{proof}
\end{proposition}

Now, we are in position to prove the main results of this section, which is a generalization of Theorem {2.7}  of [4], for $\mathcal{S}$~=~$\mathcal{S}_{I}$. Some applications of these results are appeared in Section 3.

\begin{theorem} \label{2.10}
Let $M$ be a finite $R$-module and $N$ an arbitrary $R$-module. Let $t \in \mathbb{N}_0$. Then the following conditions are equivalent:
\begin{enumerate}

\item[(i)] $Ext_R^i(M,N)\in\mathcal{S}$ for all $i\leq t$ .

\item[(ii)] For any finite $R$-module $H$ with Supp$(H)\subseteq$ Supp$(M)$,$~Ext_R^i(H,N)\in\mathcal{S}$ for all $i\leq t$.
\end{enumerate}
\begin{proof}
(i)$\Rightarrow$(ii) Since Supp$(H)\subseteq$ Supp$(M)$, according to Gruson's Theorem [{24}, Theorem 4.1], there exists a chain of submodules of $M$,
\begin{equation*}
0=H_0\subset H_1 \subset \dots \subset H_k = H,
\end{equation*} such that the factors $H_j/H_{j-1}$ are homomorphic images of a direct sum of finitely many of $M$. Now, consider the exact sequences
\begin{equation*}
0\rightarrow K \rightarrow M^n \rightarrow H_1 \rightarrow 0
\end{equation*}
\begin{equation*}
0\rightarrow H_1 \rightarrow H_2 \rightarrow H_2/H_1 \rightarrow 0
\end{equation*}
\begin{equation*}
~~~~~~~~~~~~~~~~~~\vdots
\end{equation*}
\begin{equation*}
0\rightarrow H_{k-1} \rightarrow H_k \rightarrow H_k/H_{k-1} \rightarrow 0,
\end{equation*}
for some positive integer $n$. Considering the long exact sequence
\begin{equation*}
\dots \rightarrow Ext_R^{i-1}(H_{j-1},N) \rightarrow Ext_R^i(H_j/H_{j-1},N)
\rightarrow Ext_R^i(H_j,N) \rightarrow Ext_R^i(H_{j-1},N) \rightarrow \dots
\end{equation*}
and an easy induction on $k$, the assertion follows. So, it
suffices to prove the case $k=1$. From the exact sequence
\begin{equation*}
0\rightarrow K \rightarrow M^n \rightarrow H \rightarrow 0,
\end{equation*}
where $n\in\mathbb{N}$ and $K$ is a finite $R$-module, and the
induced long exact sequence, by using induction on $i$, we show
that $Ext_R^i(H,N)\in\mathcal{S}$ for all $i$. For $i=0$, we have the exact sequence
\begin{equation*}
0\rightarrow Hom_R(H,N) \rightarrow Hom_R(M^n,N)
\rightarrow Hom_R(K,N).
\end{equation*}
Since $Hom_R(M^n,N)\cong \overset{n}{\bigoplus}
Hom_R(M,N)$, hence in view of the assumption and Lemma {2.8}, $Ext_R^0(H,N)\in\mathcal{S}$. Now, let $i>0$. We have, for any $R$-module $H$ with Supp$(H)\subseteq$ Supp$(M)$, the $R$-module
$Ext_R^{i-1}(H,N)\in\mathcal{S}$, in particular for $K$. Now, from the long exact sequence
\begin{equation*}
\dots \rightarrow Ext_R^{i-1}(K,N) \rightarrow Ext_R^i(H,N)\rightarrow
Ext_R^i(M^n,N) \rightarrow \dots
\end{equation*}
and by Lemma {2.8}, we can conclude that $Ext_R^i(H,N)\in\mathcal{S}$.\\
(ii)$\Rightarrow$(i) It is trivial.
\end{proof}
\end{theorem}

\begin{corollary} \label{2.11}
Let $r \in \mathbb{N}_0$. Then, for any $R$-module $M$, the following conditions are equivalent:
\begin{enumerate}

\item[(i)] $Ext^i_R(R/I,M)$ $\in$ $\mathcal{S}$ for all $i\leq~r$.

\item[(ii)] For any ideal $\mathfrak{a}$ of $R$ with $\mathfrak{a}\supseteq I$, $Ext^i_R(R/\mathfrak{a},M)\in\mathcal{S}$ for all $i\leq~r$.

\item[(iii)] For any finite $R$-module $N$ with Supp$(N)$ $\subseteq$ $V(I)$, $Ext^i_R(N,M)\in\mathcal{S}$ for all $i\leq~r$.

\item[(iv)] For any $\mathfrak{p}\in Min(I)$, $Ext^i_R(R/\mathfrak{p},M)\in\mathcal{S}$ for all $i\leq~r$.
\end{enumerate}
\begin{proof}
In view of Theorem {2.10}, it is enough to show that (iv) implies (i). To do this, let $\mathfrak{p}_1$,$\mathfrak{p}_2$,$\cdots$,$\mathfrak{p}_n$ be the minimal elements of $V(I)$. Then, by assumption, the $R$-modules $Ext^i_R(R/\mathfrak{p_j},M)\in\mathcal{S}$ for all $j=1,2,\cdots,n$. Hence, by Lemma \ref{2.8}, $Ext^i_R( \oplus_{j=i}^{n}R/\mathfrak{p_j},M)\cong {\oplus_{j=1}^{n}}
Ext^i_R(R/\mathfrak{p_j},M)\in\mathcal{S}$. Since Supp$(\oplus_{j=1}^{n}R/\mathfrak{p_j})$ = Supp$(R/I)$, it follows from Theorem {2.10} that $Ext^i_R(R/I,M)\in\mathcal{S}$, as required.
\end{proof}
\end{corollary}

\section{$(\mathcal{S},I,J)$-Cominimax Modules And $H_{I,J}^i(M)$}

Recall that $M$ is said to be $(I,J)$-cofinite if $M$ has support in $W(I,J)$ and $Ext^i_R(R/I,M)$ is a finite $R$-module for each $i\geq0$ (see [{23}, Definition 2.1]). In fact this definition is a generalization of $I$-cofinite modules, which is introduced by Hartshorne in [15]. Considering an arbitrary $S$erre subcategory of $R$-modules instead of finitely generated one, we can give a generalization of $(I,J)$-cofinite modules as follows.

\begin{definition} \label{3.1} Let $R$ be a Noetherian ring and $I,J$ be two ideals of $R$. For the Serre subcategory $\mathcal{S}$ of the category of $R$-modules, an $R$-module $M$ is called an $(\mathcal{S},I,J)$-cominimax precisely when Supp$(M)\subseteq W(I,J)$ and $Ext^i_R(R/I,M)\in\mathcal{S}$ for all $i\geq0$.
\end{definition}

\begin{remark} \label{3.2}
By applying various Serre classes of $R$-modules in {3.1}, we may obtain different concepts. But in view of [{22}, Proposition 1.7], the class of $(\mathcal{S},I,J)$-cominimax $R$-modules is contained in the class of $(I,J)$-torsion $R$-modules. Moreover, for every $R$-module $M$ and all $i\geq0$, $Ext^i_R(R/I,M)$ is $I$-torsion, so by Lemma {2.8} (ii), we have $Ext^i_R(R/I,M)\in\mathcal{S}_{I,J}~$ iff  $~Ext^i_R(R/I,M)\in\mathcal{S}_0$. In other words, the class of $(\mathcal{S}_0,I,J)$-cominimax $R$-modules and the class of $(\mathcal{S}_{I,J},I,J)$-cominimax $R$-modules are the same. Also, since Supp$(M)\subseteq V(I)$ implies that Supp$(M)\subseteq W(I,J)$, hence the class of $(\mathcal{S}_{I,J},I,J)$-cominimax $R$-modules contains the class of $(\mathcal{S}_I,I,J)$-cominimax $R$-modules.
\end{remark}

\begin{notation} \label{3.3}
For a Serre classes $\mathcal{S}$ of $R$-modules and two ideals $I,J$ of $R$, we use $\mathcal{C}(\mathcal{S},I,J)$ to denote the class of all $(\mathcal{S},I,J)$-cominimax $R$-modules.
\end{notation}

\begin{example and remark} \label{3.4}
(i) Let $N\in\mathcal{S}$ be such that Supp$(N)\subseteq W(I,J)$. Then it follows from Lemma {2.8} (i) that $N\in\mathcal{C}(\mathcal{S},I,J)$.

(ii) Let $N$ be a pure submodule of $R$-module $M$. By using the following exact sequence $0\rightarrow Ext_R^{i}(R/I,N) \rightarrow Ext_R^i(R/I,M)\rightarrow
Ext_R^i(R/I,M/N) \rightarrow0$, for all $i\geq0$, [{20}, Theorem 3.65], $M\in\mathcal{C}(\mathcal{S},I,J)~$ iff $~N , M/N\in\mathcal{C}(\mathcal{S},I,J)$; in particular, when $\mathcal{S}$ = $\mathcal{S}_{I,J}$ .
\end{example and remark}

\begin{proposition} \label{3.5}
let $0\rightarrow {M^\prime}\rightarrow M$ $\rightarrow M^\prime$$^\prime$ $\rightarrow0$ be an exact sequence of $R$-modules  such that two of the modules belong to $\mathcal{S}$. Then the third one is $(\mathcal{S},I,J)$-cominimax if its support is in $W(I,J)$.
\begin{proof}
The assertion follows from the induced long exact sequence
\begin{equation*}
\dots \rightarrow Ext_R^{i}(R/I,M) \rightarrow Ext_R^i(R/I,M'') \rightarrow
Ext_R^{i+1}(R/I,M^\prime)\rightarrow Ext_R^{i+1}(R/I,M)\rightarrow\dots
\end{equation*}
and Lemma {2.8} (i).
\end{proof}
\end{proposition}

An immediate consequence of Proposition {3.5} and Lemma {2.8} is as follows.

\begin{corollary} \label{3.6}
Let $f$ : $M\rightarrow N$ be a homomorphism of $R$-modules such that $M,N\in\mathcal{S}$. Let one of the three modules Ker$f$ , Im$f$ and Coker$f$ be in $\mathcal{S}$. Then two others belong to $\mathcal{C}(\mathcal{S},I,J)$ if their supports are in $W(I,J)$.
\end{corollary}

\begin{proposition} \label{3.7}
Let $I,J,I',J'$ are ideals of $R$. Then
\begin{enumerate}

\item[(i)] $M\in\mathcal{C}(\mathcal{S},I,J)$ iff $M\in\mathcal{C}(\mathcal{S},{\sqrt{I}},J)$
iff  $M\in\mathcal{C}(\mathcal{S},I,{\sqrt{J}})$ iff  $M\in\mathcal{C}(\mathcal{S},{\sqrt {I}},{\sqrt {J}})$.

\item[(ii)] If $M$ is $I$-cominimax, then  $M$ $\in$ $\mathcal{C}(\mathcal{S}_{I,J},I,J)$.

\item[(iii)] If $Min(M)$ $\subseteq W(I',J)$ , $Ext^i_R(R/I,M) \in \mathcal{S}_{I,J}$, and $I^n \subseteq {\sqrt{I'}}$ for some $n \in \mathbb{N}$ and all $i\geq0$, then $M$ $\in$ $\mathcal{C}(\mathcal{S}_{I,J},I,J)$ and $M$ $\in$ $\mathcal{C}(\mathcal{S}_{I',J},I',J)$. In particular, if $M\in\mathcal{C}(\mathcal{S}_{I,J},I,J)$, then we have $M\in\mathcal{C}(\mathcal{S}_{I',J},I',J)$.

\item[(iv)] If $Min(M)\subseteq W(I,J)$ and $J^n \subseteq {\sqrt{J'}}$ for some $n \in \mathbb{N}$, then $M\in\mathcal{C}(\mathcal{S}_{I,J},I,J)$ iff  $M\in\mathcal{C}(\mathcal{S}_{I,J'},I,J')$.
\end{enumerate}
\begin{proof}
(i) Since $V(I)=V(\sqrt I)$, the assertions follow from [{22}, Proposition 1.6], Corollary {2.11}, and Definition {3.1}.\\
(ii) By assumption and Lemma {2.8} (ii), Supp$(M)\subseteq V(I)\subseteq W(I,J)$ and $Ext^i_R(R/I,\\M)\in\mathcal{S}_{I,J}$ for all $i\geq0$.\\
(iii), (iv) Apply [{22}, Proposition 1.6 and 1.7], Corollary {2.11} and Proposition {2.7}(iii),\\(iv).
\end{proof}
\end{proposition}

The following Remark plays an important role in the proof of our main theorems in this section.

\begin{remark} \label{3.8}
In view of proof [{22}, Theorem 3.2], $\Gamma_\mathfrak{a}(M)\subseteq\Gamma_{I,J}(M)$, for any $\mathfrak{a}\in\tilde{W}(I,J)$. Thus $\Gamma_{I,J}(M)$= 0 implies that $\Gamma_\mathfrak{a}(M)$ = 0, for all $\mathfrak{a}\in\tilde{W}(I,J)$. Now, let $\bar{M}$ = $M/\Gamma_{I,J}(M)$ and $E$ = $E_{R}(\bar{M})$ be the injective hull of $R$-module $\bar{M}$. Put $L$ = $E/\bar{M}$. Since $\Gamma_{I,J}(\bar{M})$ = 0, then $\Gamma_{I,J}(E)$ = 0 and also for any $\mathfrak{a}\in\tilde{W}(I,J)$, we have $\Gamma_\mathfrak{a}(\bar{M})$ = 0 = $\Gamma_\mathfrak{a}(E)$. In particular, the $R$-module $Hom_R(R/\mathfrak{a},E)$ is zero. Now, from the exact sequence $0\rightarrow \bar{M}\rightarrow E \rightarrow L \rightarrow 0$, and applying $Hom_R(R/\mathfrak{a},-)$ and $\Gamma_{I,J}(-)$, we have the following isomorphisms
$$ Ext^i_R(R/\mathfrak{a},L)\cong Ext^{i+1}_R(R/\mathfrak{a},\bar{M}) ~and ~H^i_{I,J}(L)\cong H^{i+1}_{I,J}(M), $$
for any $\mathfrak{a}\in\tilde{W}(I,J)$ and all $i\geq0$. In particular, $Ext^i_R(R/I,L)\cong Ext^{i+1}_R(R/I,\bar{M})$.
\end{remark}

\begin{proposition} \label{3.9}
Let $t\in\mathbb{N}_0$ be such that $H^i_{I,J}(M)\in\mathcal{C}(\mathcal{S},I,J)$ for all $i<t$. Then $Ext^i_R(R/I,M)\in\mathcal{S}$ for all $i<t$.
\begin{proof}
We use induction on $t$. When $t$
= 0, there is nothing to prove. For $t$ = 1, since $Hom_R(R/I,\Gamma_{I,J}(M))$ = $Hom_R(R/I,M)$, and $\Gamma_{I,J}(M)$ is $(\mathcal{S},I,J)$-cominimax, the result is true. Now, suppose that $t\geq2$ and the case $t-1$ is settled. The exact sequence $0\rightarrow\Gamma_{I,J}(M)\rightarrow M\rightarrow
\bar{M}\rightarrow 0$ induced the long exact sequence
\begin{equation*}
\dots \rightarrow Ext_R^{i}(R/I,\Gamma_{I,J}(M)) \rightarrow Ext_R^i(R/I,M) \rightarrow
Ext_R^{i}(R/I,\bar{M}) \rightarrow \dots~.
\end{equation*}
Since $\Gamma_{I,J}(M)$ $\in$ $\mathcal{C}(\mathcal{S},I,J)$, we have $Ext^i_R(R/I,\Gamma_{I,J}(M))\in\mathcal{S}$ for all $i\geq0$. Therefore, it is enough to show that $Ext^i_R(R/I,\bar{M})\in\mathcal{S}$ for all $i<t$. For this purpose, let $E$ = $E_R(\bar{M})$ and $L$ = $E/\bar{M}$. Now, by Remark {3.8}, for all $i\geq0$, we get the isomorphisms $H^i_{I,J}(L)\cong H^{i+1}_{I,J}(M)$ and  $Ext^i_R(R/I,L)\cong Ext^{i+1}_R(R/I,\bar{M})$. Now, by assumption, $H^{i+1}_{I,J}(M)\in\mathcal{C}(\mathcal{S},I,J)$ for all $i<t-1$, and so $H^i_{I,J}(L)\in\mathcal{C}(\mathcal{S},I,J)$ for all $i<t-1$. Thus, by the inductive hypothesis, $Ext^i_R(R/I,L)\in\mathcal{S}$ and so $Ext^{i+1}_R(R/I,\bar{M})\in\mathcal{S}$.
\end{proof}
\end{proposition}

The next corollary generalizes Proposition {3.7} of [4].

\begin{corollary} \label{3.10}
Let $H^i_{I,J}(M)\in\mathcal{C}(\mathcal{S},I,J)$ for all $i\geq0$. Then $Ext^i_R(R/I,M)\in\mathcal{S}$ for all $i$ $\geq$$0$; particularly, when $\mathcal{S}$ is the class of $I$-minimax modules or the class of $(I,J)$-minimax modules.
\end{corollary}

The Proposition {3.8} of [4] can be obtained from the following theorem when $J=0$ and $\mathcal{S}=\mathcal{S}_{I}$.

\begin{theorem} \label{3.11}
Let $Ext^i_R(R/I,M)\in\mathcal{S}$ for all $i\geq0$. Let $t\in\mathbb{N}_0$ be such that $H^i_{I,J}(M)\in\mathcal{C}(\mathcal{S},I,J)$, for all $i\neq t$, then $H^t_{I,J}(M)\in\mathcal{C}(\mathcal{S},I,J)$.
\begin{proof}
We use induction on $t$. If $t=0$, we must prove that $Ext^i_R(R/I,\Gamma_{I,J}(M))\in\mathcal{S}$ for all $i\geq0$. By the exact sequence \begin{equation*}
\dots \rightarrow Ext_R^{i-1}(R/I,\bar{M}) \rightarrow Ext_R^i(R/I,\Gamma_{I,J}(M)) \rightarrow
Ext_R^{i}(R/I,M) \rightarrow \dots
\end{equation*}
and the hypothesis, it is enough to show that $Ext^i_R(R/I,\bar{M})\in\mathcal{S}$ for all $i\geq0$. Now, by Remark {3.8} and our assumption, we obtain $H^i_{I,J}(L)\in\mathcal{C}(\mathcal{S},I,J)$. Therefore Corollary {3.10} implies that $Ext^i_R(R/I,\bar{M})\in\mathcal{S}$ for all $i\geq0$ (note that $Ext^0_R(R/I,\bar{M})=0)$. Now suppose, inductively, that $t>0$ and the result has been proved for $t-1$. By Remark {3.8}, it is easy to show that $L$ satisfies in our inductive hypothesis. Therefore, the assertion follows from $H^t_{I,J}(M)$
$\cong$ $H^{t-1}_{I,J}(L)$.
\end{proof}
\end{theorem}

\begin{corollary} \label{3.12}
Let $M\in\mathcal{S}$ and  $t\in\mathbb{N}_0$ be such that $H^i_{I,J}(M)$ is $(\mathcal{S},I,J)$-cominimax for all $i\neq t$.
Then $H^t_{I,J}(M)$ is $(\mathcal{S},I,J)$-cominimax.
\begin{proof}
This is an immediate consequence of Lemma {2.8} (i) and Theorem {3.11}.
\end{proof}
\end{corollary}

\begin{corollary} \label{3.13}
Let $I$ be a principal ideal and $J$ be an arbitrary ideal of $R$. Let $M$ $\in$ $\mathcal{S}$. Then $H^i_{I,J}(M)$ is $(\mathcal{S},I,J)$-cominimax for all $i\geq0$.
\begin{proof}
For $i=0$, since $H^0_{I,J}(M)$ is a submodule of $M$ and $M\in\mathcal{S}$, it turns out that $H^0_{I,J}(M)$ is $(\mathcal{S},I,J)$-cominimax, by Remark {3.4} (i). Now, let $I=aR$. By [{22}, Definition 2.2 and Theorem 2.4], we have $H^i_{I,J}(M)\cong H^i(C^\bullet_{I,J}\otimes_RM)=0$ for all $i>1$. Therefore the result follows from Theorem {3.11}.
\end{proof}
\end{corollary}

Now we are prepared to prove the main theorem of this section, which is a generalization of one of the main results of [{3}, Theorem 2.2] and also [{23}, Theorem 2.3].

\begin{theorem} \label{3.14}
Let  $\mathfrak{a}\in\tilde{W}(I,J)$. Let $t\in\mathbb{N}_0$ be such that $Ext^t_R(R/\mathfrak{a},M)\in\mathcal{S}$ and $Ext^j_R(R/\mathfrak{a},H^i_{I,J}(M))\in\mathcal{S}$ for all $i<t$ and all $j\geq0$. Then for any submodule $N$ of $H^t_{I,J}(M)$ such that $Ext^1_R(R/\mathfrak{a},N)\in\mathcal{S}$, we have $Hom_R(R/\mathfrak{a},H^t_{I,J}(M)/N)$
$\in$ $\mathcal{S}$; in particular, for $\mathfrak{a}=I$.
\begin{proof}
Considering the following long exact sequence \begin{equation*}
\dots \rightarrow Hom_R(R/\mathfrak{a},H^t_{I,J}(M)) \rightarrow Hom_R(R/\mathfrak{a},H^t_{I,J}(M)/N) \rightarrow
Ext_R^1(R/\mathfrak{a},N) \rightarrow \dots,
\end{equation*}
since $Ext_R^1(R/\mathfrak{a},N)\in\mathcal{S}$, it is enough to show that $Hom_R(R/\mathfrak{a},H^t_{I,J}(M))\in\mathcal{S}$. To do this, we use induction on $t$. When $t=0$, since $Hom_R(R/\mathfrak{a},\Gamma_{I,J}(M))$ = $Hom_R(R/\mathfrak{a},M)\in\mathcal{S}$, the result is obtained. Next, we assume that $t>0$ and that the claim is true for $t-1$. Let $\bar{M}$ = $M/\Gamma_{I,J}(M)$. Then, by the long exact sequence \begin{equation*}
\dots \rightarrow Ext_R^j(R/\mathfrak{a},M) \rightarrow Ext_R^j(R/\mathfrak{a},\bar{M})\rightarrow
Ext_R^{j+1}(R/\mathfrak{a},\Gamma_{I,J}(M)) \rightarrow \dots,
\end{equation*}
and assumption, we conclude that $Ext^j_R(R/\mathfrak{a},\bar{M})\in\mathcal{S}$. Now, by using notation of Remark {3.8}, it is easy to see that $L$ satisfies the inductive hypothesis. So that we get $Hom_R(R/\mathfrak{a},H^{t-1}_{I,J}(L))\in\mathcal{S}$ and therefore,  $Hom_R(R/\mathfrak{a},H^t_{I,J}(M))\in\mathcal{S}$, as required.
\end{proof}
\end{theorem}

The main results of [{4}, Theorem 4.2] , [{5}, Lemma 2.2] , [{2}, Theorem 1.2] ,[16] , [{12}, Corollary 2.7], and [{17}, Corollary 2.3] are all special cases of next corollary, by replacing various Serre classes with $\mathcal{S}$ and $J$ = 0.

\begin{corollary} \label{3.15}
Let $t\in\mathbb{N}_0$ be such that $Ext^t_R(R/I,M)\in\mathcal{S}$ and $H^i_{I,J}(M)\in\mathcal{C}(\mathcal{S},I,J)$ for all $i<t$. Then for any submodule $N$ of $H^t_{I,J}(M)$ and any finite $R$-module $M'$ with $Supp(M')\subseteq V(I)$ and $Ext^1_R(M',N)\in\mathcal{S}$, we have $Hom_R(M',\\H^t_{I,J}(M)/N)\in\mathcal{S}$.
\begin{proof}
Apply Theorem {3.14} and Corollary {2.11}.
\end{proof}
\end{corollary}

\begin{proposition} \label{3.16}
Let $t\in\mathbb{N}_0$ be such that $H^i_{I,J}(M)$ $\in$ $\mathcal{C}(\mathcal{S},I,J)$ for all $i<t$. Then the following statements hold:

\begin{enumerate}

\item[(i)] If $Ext^t_R(R/I,M)\in\mathcal{S}$, then $Hom_R(R/I,H^t_{I,J}(M))\in\mathcal{S}$.

\item[(ii)] If $Ext^{t+1}_R(R/I,M)\in\mathcal{S}$, then $Ext^1_R(R/I,H^t_{I,J}(M))\in\mathcal{S}$.

\item[(iii)] If $Ext^i_R(R/I,M)\in\mathcal{S}$ for all $i\geq0$, then $Hom_R(R/I,H^{t+1}_{I,J}(M))\in\mathcal{S}$ iff $Ext^2_R(R/I,H^t_{I,J}(M))\in\mathcal{S}$.
\end{enumerate}

\begin{proof}
(i) Apply Corollary {3.15} or Theorem {3.14} .\\
(ii) We proceed by induction on $t$. If $t=0$, then by the long exact sequence
\begin{enumerate}
\item[]$(\ast)$ $~~$ $~~$ $~0\rightarrow Ext_R^1(R/I,\Gamma_{I,J}(M)) \rightarrow Ext_R^1(R/I,M)\rightarrow Ext_R^1(R/I,\bar{M})$\\

\item[]  $~~$ $~~$ $~~$ $~~$ $~~$ $~~$ $\rightarrow Ext_R^2(R/I,\Gamma_{I,J}(M)) \rightarrow Ext_R^2(R/I,M)\rightarrow Ext_R^2(R/I,\bar{M})$

\item[]  $~~$ $~~$ $~~$ $~~$ $~~$ $~~$ $\vdots$

\item[]  $~~$ $~~$ $~~$ $~~$ $~~$ $~~$ $\rightarrow Ext_R^i(R/I,\Gamma_{I,J}(M)) \rightarrow Ext_R^i(R/I,M)\rightarrow Ext_R^i(R/I,\bar{M})$

\item[]  $~~$ $~~$ $~~$ $~~$ $~~$ $~~$ $\rightarrow Ext_R^{i+1}(R/I,\Gamma_{I,J}(M))\rightarrow \cdots$\\
\end{enumerate}
and $Ext^1_R(R/I,\bar{M})\in\mathcal{S}$, the result follows. Suppose that $t>0$ and the assertion is true for $t-1$. Since $\Gamma_{I,J}(M)\in\mathcal{C}(\mathcal{S},I,J)$, so $Ext^i_R(R/I,\Gamma_{I,J}(M))\in\mathcal{S}$ for all $i\geq0$, and so by $(\ast)$, $Ext^{t+1}_R(R/I,\bar{M})\in\mathcal{S}$. Now, by the notations of Remark {3.8}, it is easy to see that $R$-module $L$ satisfies the inductive hypothesis and so $Ext^1_R(R/I,H^{t-1}_{I,J}(M))\in\mathcal{S}$. Now, the result follows from $H^{t-1}_{I,J}(L)$ $\cong$ $H^{t}_{I,J}(M)$.\\
(iii) $(\Rightarrow)$ We use induction on $t$. Let $t=0$. Then considering the long exact sequence $(\ast)$, it is enough to show that $Ext^1_R(R/I,\bar{M})\in\mathcal{S}$. By Remark {3.8}, we have

$$Ext_R^1(R/I,\bar{M})\cong Hom_R(R/I,L)$$
$$\ \ \ \ \ \ \ \ \ \ \ \ \ \ \ \ \ \ \ \ \ \ \ \ \ \ \ \ \cong Hom_R(R/I,\Gamma_{I,J}(L))$$
$$ \ \ \ \ \ \ \ \ \ \ \ \ \ \ \ \ \ \ \ \ \ \ \ \ \ \ \ \ \ \ \cong Hom_R(R/I,H^1_{I,J}(M)),$$
as required. Suppose $t>0$ and the assertion is true for $t-1$. Since $\Gamma_{I,J}(M)$ $\in$ $\mathcal{C}(\mathcal{S},I,J)$, we have $Ext^i_R(R/I,\Gamma_{I,J}(M))$ $\in$ $\mathcal{S}$ for all $i\geq0$. Therefore the exactness of sequence $(\ast)$ implies that $Ext^i_R(R/I,\bar{M})\in\mathcal{S}$ for all $i\geq0$. Again by using the notations of Remark {3.8}, we get $Ext^i_R(R/I,L)$ $\in$ $\mathcal{S}$, for all $i\geq0$, and also  $Hom_R(R/I,H^t_{I,J}(L))\cong Hom_R(R/I,H^{t+1}_{I,J}(M))$ $\in$ $\mathcal{S}$. Now, by inductive hypothesis, $Ext^2_R(R/I,H^{t-1}_{I,J}(L))\in\mathcal{S}$ and hence $Ext^2_R(R/I,H^t_{I,J}(M))\in\mathcal{S}$, as required.\\
$(\Leftarrow)$ This part can be proved by the same method of $(\Rightarrow)$, using induction on $t$, the following exact sequence
\begin{equation*}
Ext_R^1(R/I,M)\rightarrow Ext_R^1(R/I,\bar{M})\rightarrow Ext_R^2(R/I,\Gamma_{I,J}(M)),
\end{equation*} and Remark {3.8}.
\end{proof}
\end{proposition}

\section{Finiteness Properties Of Associated  Primes}

In this short section, we obtain some results, as some applications of previous sections, about associated prime ideals of local cohomology modules and also finiteness properties of them.

\begin{proposition} \label {4.1}
Let $t \in \mathbb{N}_0$ be such that $Ext^t_R(R/I,M)$ $\in$ $\mathcal{S}_{I,J}$ and $H^i_{I,J}(M)$ $\in$ $\mathcal{C}(\mathcal{S}_{I,J},I,J)$ for all $i<t$. Let $N$ be a submodule of $H^t_{I,J}(M)$ such that $Ext^1_R(R/I,N)$ $\in$ $\mathcal{S}_{I,J}$. If Supp$(H^t_{I,J}(M)/N)\subseteq V(I)$, then Gdim$(H^t_{I,J}(M)/N)<\infty$ and so $H^t_{I,J}(M)/N$ has finitely many associated primes.
\begin{proof}
By using Theorem {3.14}, for the Serre class $\mathcal{S}_{I,J}$, we have $Hom_R(R/I,\\H^t_{I,J}(M)/N)$ $\in$ $\mathcal{S}_{I,J}$. Hence, by Lemma {2.8} (ii), $Hom_R(R/I,H^t_{I,J}(M)/N)$ $\in$ $\mathcal{S}_0$, as required.
\end{proof}
\end{proposition}

\begin{corollary} \label {4.2}
Let $t \in \mathbb{N}_0$ be such that Ext$^t_R(R/I,M)$ and $H^i_{I,J}(M)$ are $(I,J)$-minimax $R$-modules for all $i<t$. Let $N$ be a submodule of $H^t_{I,J}(M)$ such that Supp$(H^t_{I,J}(M)/N)\subseteq V(I)$ and Ext$^1_R(R/I,N)$ is $(I,J)$-minimax. Then $H^t_{I,J}(M)/N$ has finite Goldie dimension and so Ass$(H^t_{I,J}(M)/N)$ is a finite set; in particular for $N$ = $JH^t_{I,J}(M)$.
\begin{proof}
For the first part, apply Remark {3.4} and proposition {4.1}. Since by [{22}, Corollary 1.9], $H^i_{I,J}(M)/JH^t_{I,J}(N)$ is $I$-torsion, so the last part immediately follows from the first.
\end{proof}
\end{corollary}

\begin{corollary} \label {4.3}
(See [{21}, Theorem 4]) Let $t \in \mathbb{N}_0$ be such that $Ext^t_R(R/I,M)$ is a finite $R$-module. If $H^i_I(M)$ is $I$-cofinite for all $i<t$ and $H^t_I(M)$ is minimax, then $H^t_I(M)$ is $I$-cofinite and so Ass$((H^t_I(M))$ is a finite set.
\begin{proof}
In Proposition {3.16} (i), apply $\mathcal{S}$ as the class of finite $R$-modules, and $J$ = 0. Therefore, $Hom_R(R/I,H^t_I(M))$ is finite $R$-module. Now, use [{18}, Proposition 3.4].
\end{proof}
\end{corollary}

\begin{corollary} \label {4.4}
Let the situation be as in Corollary 4.3. Then the following statements hold:

\begin{enumerate}

\item[(i)] If $Ext^{t+1}_R(R/I,M)$ is finite, then $Hom_R(R/I,H^{t+1}_I(M))$ and  $Ext^1_R(R/I,H^t_I(M))$ are finite and so Ass$(H^{t+1}_I(M))$ is a finite set.

\item[(ii)] If $Ext^i_R(R/I,M)$ is finite for all $i\geq0$, then  $Ext^2_R(R/I,H^t_{I,J}(M))$ is finite.

\end{enumerate}

\begin{proof}
(i) By Corollary 4.3, we conclude that $H^t_I(M)$ is $I$-cofinite. So $H^i_I(M)$ is $I$-cofinite for all $i<t+1$. Now, using Proposition {3.16} (i), (ii).\\
(ii) The result follows from (i) and Proposition {3.16} (iii).
\end{proof}
\end{corollary}


\bibliographystyle{amsplain}

\end{document}